\input amstex
\documentstyle{amsppt}
\NoBlackBoxes \magnification=\magstep1 \vsize=22 true cm \hsize=16.1
true cm \voffset=1 true cm
%%%%%%%%%%%%%%%%%%%%%%%%%%%%%%%%%%%%%%%%%%%%%%%%%%%%%%%%%%%%

\redefine\id{\operatorname{id}}
\redefine\BMO{\operatorname{BMO}}
\redefine\VMO{\operatorname{VMO}}
\redefine\WP{\operatorname{WP}}

\redefine\Diff{\operatorname{Diff}}
\redefine\Rot{\operatorname{Rot}}
\topmatter
\title     Weil-Petersson Teichm\"{u}ller space
\endtitle
\title     Weil-Petersson Teichm\"{u}ller space II: smoothness of flow curves of $H^{\frac 32}$-vector fields
\endtitle
\author    Yuliang Shen  \quad  Shuan Tang
\endauthor
%\affil (Department of Mathematics,  Soochow  University)
%\endaffil
\address
$$
\align
&{\text{Yuliang Shen}}\\
&{\text{Department of Mathematics,  Soochow  University,   Suzhou
215006,  P.\,R. China}}\\
&{\text{email: ylshen\@suda.edu.cn}}\\
&\\
&{\text{Shuan Tang}}\\
&{\text{School of Mathematical Sciences, Guizhou Normal University, Guiyang 550001, P.\,R. China}}\\
&{\text{email: tsa\@gznu.edu.cn}}\\
\endalign
$$
\endaddress
\abstract Given  a continuous vector field $\lambda(t, \cdot)$ of Sobolev class $H^{\frac 32}$ on the unit circle $S^1$,  the flow maps $\eta=g(t, \cdot)$ of the differential equation
$$
\cases
\frac{d\eta}{dt}=\lambda(t, \eta)\\
\eta(0,\zeta)=\zeta
\endcases
$$
are known to be quasisymmetric homeomorphisms. Very recently, Gay-Balmaz-Ratiu [GR] conjectured that the flow curve $g(t, \cdot)$ is in the Weil-Petersson class WP$(S^1)$ and is continuously differentiable with respect to the Hilbert manifold structure of WP$(S^1)$ introduced by Takhtajan-Teo [TT]. The first assertion had already been  demonstrated in our previous paper [Sh2]. In this sequel to [Sh2], we will continue to deal with the Weil-Petersson class WP$(S^1)$ and completely solve this conjecture in the affirmative.
\endabstract
\thanks
{\it 2010 Mathematics Subject Classification}: Primary 30C62; 30F60; 32G15, Secondary 30H30; 30H35; 46E35
\endthanks
\thanks
{\it Key words and phrases}: Universal Teichm\"uller space; Weil-Petersson
Teichm\"uller space; quasi-symmetric homeomorphism;
quasiconformal mapping;   Sobolev class
\endthanks
\thanks Research  supported by  the National Natural Science Foundation of China (Grant Nos. 11631010, 11601100).
\endthanks
\endtopmatter

\document

\baselineskip 13 pt

\head 1 Introduction\endhead

This  is a continuous work of our previous paper [Sh2], where we  presented some recent results on the Weil-Petersson geometry  theory of the universal Teichm\"uller space, a topic which  is important in Teichm\"uller theory and has wide applications to various areas such as mathematical physics (see [BR1-2], [Ki], [KY], [RSW1-2]),  differential equation and computer vision (see [GMR], [GR], [Ku]).

A sense-preserving homeomorphism $g$  of the unit circle $S^1$ onto itself is said to belong to  the Weil-Petersson class, which is denoted by $\WP(S^1)$, if it  has a quasiconformal extension  to
the unit disk $\Delta$ whose Beltrami coefficient $\nu$ is square integrable in the Poincar\'e metric, namely,
$$\iint_{\Delta}|\nu(z)|^2(1-|z|^2)^{-2}dxdy<+\infty.\tag 1.1$$
 In an important paper [TT], Takhtajan-Teo  showed how to endow $\WP(S^1)$ and its two close relatives, the Weil-Petersson Teichm\"uller space $T_0=\WP(S^1)/\text{M\"ob}(S^1)$ and the  Weil-Petersson Teichm\"uller curve $\hat T_0=\WP(S^1)/\Rot(S^1)$, with Hilbert manifold structures (see also [GR], [Sh2]). Here, M\"ob$(S^1)$ denotes the group of all M\"obius transformations keeping the unit disk $\Delta$ fixed, while $\Rot(S^1)$ denotes the sub-group of all rotations about the circle $S^1$. In [Sh2], we gave the following intrinsic characterization  of a quasisymmetric homeomorphism in the Weil-Petersson class $\WP(S^1)$ without using quasiconformal extensions, which solves a problem proposed by Takhtajan-Teo in 2006. Recall that, for a function $f$ defined on a subset $\Gamma$ of the complex plane, $f'$ denotes the derivative of  $f$, namely, for $z\in \Gamma$,
$$f'(z)\doteq\lim_{\Gamma\ni \zeta\to z}\frac{f(\zeta)-f(z)}{\zeta-z} $$
provided the limit exists, while $f'(z)\doteq 0$ otherwise.

 \proclaim{Theorem 1.1 ([Sh2])} A sense-preserving  homeomorphism $g$ on the unit circle $S^1$ belongs to the Weil-Petersson class $\WP(S^1)$  if and only if $g$ is absolutely continuous  $($with respect to the arc-length measure$)$ such that $\log g'$ belongs to the Sobolev  class $H^{\frac 12}$ on the unit circle. Moreover, the correspondence $g\mapsto\log |g'|$ induces a homeomorphism from the  Weil-Petersson Teichm\"uller curve $\hat T_0$ onto the real Sobolev space $H_{\Bbb R}^{\frac 12}/\Bbb R$.
 \endproclaim

 It should be pointed out that the second assertion  of Theorem 1.1 was stated in a different but an equivalent way in [Sh2] (see Theorem 8.1 in [Sh2]). The definition of Sobolev class  will be given in the next section. In this paper, we will continue to deal with the Weil-Petersson class WP$(S^1)$. Recall that WP$(S^1)$ is modeled on the Sobolev  space $H^{\frac 32}$, namely, the tangent space to WP$(S^1)$ at the identity consists of precisely the $H^{\frac 32}$ vector fields $\lambda$ on the unit circle (see [NV], [TT] and also [HSWS]). We will be mainly concerned with  the  flows of $H^{\frac 32}$ vector fields on the unit circle.
 It is known that the Weil-Petersson class $\WP(S^1)$ can be generated by the flows of the $H^{\frac 32}$ vector fields on the unit circle (see [Fi], [GR]).
Here we consider the converse problem and  prove the following result, completely solving a  conjecture posed by  Gay-Balmaz-Ratiu in the recent paper [GR] (page 760).

\proclaim{Theorem 1.2} Let $\lambda(t, \cdot)\in C^0([0, M], H^{\frac 32})$ be  a continuous vector field  of Sobolev class $H^{\frac 32}$ on the unit circle $S^1$. Then  the flow curve $\eta=g(t, \cdot)$ of the differential equation
$$
\cases
\frac{d\eta}{dt}=\lambda(t, \eta)\\
\eta(0,\zeta)=\zeta
\endcases\tag 1.2
$$
 is in the Weil-Petersson class $\WP(S^1)$ and is continuously differentiable with respect to the Hilbert manifold structure of $\WP(S^1)$ such that
 $$\frac{d}{dt}g(t, \cdot)=\lambda(t, g(t, \cdot)).\tag 1.3$$
 \endproclaim

Recall that the first assertion in Theorem 1.2 was already proved by the author in [Sh2]. In fact, we have proved

 \proclaim{Theorem 1.3 ([Sh2])} Under the assumption of Theorem 1.2, the flow curve $g(t, \cdot)$ of the differential equation (1.2) satisfies $\log g'(t, \cdot)\in H^{\frac 12}$, which implies by Theorem 1.1 that the flow curve $\eta=g(t, \cdot)$ is in the Weil-Petersson class $\WP(S^1)$,  and  the mapping $t\mapsto \log g'(t, \cdot)$ from $[0, M]$ into $H^{\frac 12}$ is continuously differentiable such that
$$\frac{d}{dt}\log g'(t, \cdot)=\lambda'(t, g(t, \cdot)).\tag 1.4$$
\endproclaim

Theorem 1.2 has several important consequences on the regularity of the Weil-Petersson class $\WP(S^1)$ and on the flows of the vector fields  of Sobolev class $H^{\frac 32}$ on the unit circle $S^1$ (see [GR] and [Sh2] for more details). It  is also assumed  to be useful to the further study of the geometry  of the Weil-Petersson Teichm\"uller space $T_0$. We hope to persue this in a separated paper.

An open problem  (see page 68 in [TT]) is to give a geometric characterization of a Weil-Petersson quasi-circle, the image of the unit circle $S^1$ under a quasiconformal mapping which is conformal outside the unit disk $\Delta$ and has Beltrami coefficient in $\Delta$ satisfying (1.1). A partial answer to this problem  was obtained by Gallardo-Guti\'errez,  Gonz\'alez,  P\'erez-Gonz\'alez,  Pommerenke and  R\"atty\"a [GGPPR]. A Weil-Petersson quasi-line is defined in the same way, namely, it is the image of the real line $\Bbb R$  under a quasiconformal mapping which is conformal on the lower plane $\Bbb U^*$  and has Beltrami coefficient in the upper half plane $\Bbb U$ being square integrable in the Poincar\'e metric, that is, satisfying (2.1) below. In a forth-coming paper [SW], we will endow the set of all Weil-Petersson quasi-lines (with certain normalized condition) with a real Hilbert manifold structure from a geometric point of view and show that that new manifold structure is topologically equivalent to the standard complex Hilbert manifold structure given by Takhtajan-Teo [TT]. Theorem 2.3 in the next section  will play an essential role in that work.

 \head 2 Weil-Petersson Teichm\"uller space on the real line
\endhead

In this section, we give some basic definitions and results on the Weil-Petersson Teichm\"uller space (see [Sh2] and [TT] for more details). As will be seen later, it is convenient to define the Weil-Petersson Teichm\"uller space and first prove Theorem 1.2 in the setting of the real line $\Bbb R$ instead of the unit circle $S^1$. Actually, as stated at the end of the first section, the results in the real line case, eg. Theorems 2.2 and  2.3 below,   turn out to be very useful to the study of geometric characterizations of Weil-Petersson quasi-lines (see [SW]).

Let  $M(\Bbb U)$ denote  the open
unit ball of the Banach space $L{}^{\infty}(\Bbb U)$ of
essentially bounded measurable functions on the upper half plane $\Bbb U$ in the complex plane $\Bbb C$. For $\mu\in
M(\Bbb U)$, let $f_{\mu}$ be the unique quasiconformal mapping on $\Bbb U$ onto itself which has complex dilatation $\mu$ and keeps the points $0$, $1$ and $\infty$ fixed.  We say two elements $\mu$
and $\nu$ in $M(\Bbb U)$ are equivalent, denoted by $\mu\sim\nu$,
if $f_{\mu}=f_{\nu}$ on the real line $\Bbb R$. Then
$T=M(\Bbb U)/_{\sim}$ is known as the Bers model of the universal
Teichm\"uller space. We let $\Phi$ denote the natural  projection
 from $M(\Bbb U)$ onto $T$ so that $\Phi(\mu)$ is the equivalence class $[\mu]$. $[0]$ is called the base point of $T$.  It is known that ${T}$ has a unique   complex Banach manifold structure such that
 the
natural projection $\Phi$ from $M(\Bbb U)$ onto ${T}$ is a
holomorphic split submersion (see [GL], [Le], [Na1]).

  We denote by
$\Cal L(\Bbb U)$ the Banach space of all
measurable functions $\mu$ with norm
$$\|\mu\|_{\WP}\doteq\|\mu\|_{\infty}+\left(\frac{1}{\pi}\iint_{\Bbb U}\frac{|\mu(z)|^2}{y^2}dxdy\right)^{\frac 12}, \quad z=x+iy.\tag 2.1$$
 Set $\Cal
M(\Bbb U)=M(\Bbb U)\cap\Cal L(\Bbb U)$. Then ${T_0}=\Cal
M(\Bbb U)/_{\sim}$ is the complex model of the Weil-Petersson Teichm\"uller space. It is known that  ${T_0}$ has a unique   complex Hilbert manifold structure such that
 the
natural projection $\Phi$ from $\Cal M(\Bbb U)$ onto ${T_0}$ is a
holomorphic split submersion (see [TT] and also [Sh2]).

It is well known that a quasiconformal self-mapping of $\Bbb U$
 induces a bi-holomorphic automorphism of the universal
Teichm\"uller space (see [GL], [Le], [Na1]). Precisely,  let
$w:\Bbb U\to\Bbb U$ be a quasiconformal mapping with complex
dilatation $\mu$. Then $w$ induces an bi-holomorphic isomorphism
$R_w: M(\Bbb U)\to M(\Bbb U)$ as
$$R_w(\nu)=\left(\frac{\nu-\mu}{1-\bar\mu\nu}\frac{\partial w}{\overline{\partial w}}\right)\circ w^{-1}.\tag 2.2$$ $R_w$
descends down a bi-holomorphic isomorphism $w^*: T\to T$ by
$w^*\circ\Phi=\Phi\circ R_w$. $w^*$ is usually called an allowable mapping.
When $w$ is quasi-isometric under the
Poincar\'e metric $|dz|/y$ with Beltrami coefficient  $\mu\in\Cal
M(\Bbb U)$, $R_w$ maps $\Cal M(\Bbb U)$ into itself and $w^*: {T_0}\to {T_0}$ is bi-holomorphic.

We denote by $\WP(\Bbb R)$ the Weil-Petersson class of all increasing  homeomorphisms $h$ of $\Bbb R$ onto itself which have quasiconformal extensions $w$ to the upper half plane $\Bbb U$ whose Beltrami coefficients $\mu$ belong to the class $\Cal M(\Bbb U)$. We also denote by  $\WP_0(\Bbb R)$ the sub-class of $\WP(\Bbb R)$ of all mappings $h$ with the normalized condition $h(0)=0$, $h(1)=1$. Then the correspondence  $[\mu]\mapsto f_{\mu}|_{\Bbb R}$ induces a one-to-one map $I$ from $T_0$ onto the normalized Weil-Petersson class $\WP_0(\Bbb R)$, which endows $\WP_0(\Bbb R)$  with a complex Hilbert manifold structure (under which $I$ is a bi-holomorphic isomorphism).

Recall that the Sobolev class $H^{\frac 12}$ ($H_{\Bbb R}^{\frac 12}$)  on the unit circle $S^1$ or the real line $\Bbb R$ is the set of all locally integrable (real-valued) functions $\varphi$ with
$$\|\varphi\|^2_{H^{\frac 12}}=\frac{1}{4\pi^2}\int_{\Bbb S}\int_{\Bbb S}\frac{|\varphi(u)-\varphi(\tilde u)|^2}{|u-\tilde u|^2}|du||d\tilde u|<+\infty,\tag 2.3$$
where $\Bbb S$ denotes the unit circle $S^1$ or the real line $\Bbb R$.
We denote by $H^{\frac 32}$ ($H_{\Bbb R}^{\frac 32}$)  the class of all (real-valued) functions $\varphi$ on the  unit circle $S^1$ or the real line $\Bbb R$ which are  locally absolutely continuous such that $\varphi'\in H^{\frac 12}$ ($H_{\Bbb R}^{\frac 12}$). As will be seen in section 8 (Theorem 8.1 below), the tangent space to $\WP_0(\Bbb R)$ at the identity consists of precisely the $H^{\frac 32}$ real-valued  vector fields  on the real line vanishing at the points $0$ and $1$.

We have the following result on the real line parallel to Theorem 1.2.

\proclaim{Theorem 2.1} Let $\omega(t, \cdot)\in C^0([0, M], H_{\Bbb R}^{\frac 32})$ be  a continuous real-valued vector field   on the real line $
\Bbb R$ with the normalized condition $\omega(t, 0)=\omega(t, 1)=0$. Then  the flow curve $u=h(t, \cdot)$ of the differential equation
$$
\cases
\frac{du}{dt}=\omega(t, u)\\
u(0,x)=x
\endcases\tag 2.4
$$
 is in the normalized Weil-Petersson class $\WP_0(\Bbb R)$ and is continuously differentiable with respect to the Hilbert manifold structure of $\WP_0(\Bbb R)$ such that
 $$\frac{d}{dt}h(t, \cdot)=\omega(t, h(t, \cdot)).\tag 2.5$$
 \endproclaim

The following result plays an essential role in the proof of Theorem 2.1.

\proclaim{Theorem 2.2} Let $h$ be an increasing    and  locally absolutely continuous homeomorphism from  the real line  onto itself  such that $\log h'$ belongs to the Sobolev  class $H^{\frac 12}$. Then $h$  belongs to the Weil-Petersson class $\WP(\Bbb R)$. Moreover, the correspondence
 $$u\mapsto h_u:\quad h_u(x)=\frac{\int_0^xe^{u(t)}dt}{\int_0^1e^{u(t)}dt},\quad x\in\Bbb R\tag 2.6$$
 induces a real analytic map $\Psi$  from the  real Sobolev space $H_{\Bbb R}^{\frac 12}/\Bbb R$ into the normalized Weil-Petersson class $\WP_0(\Bbb R)(=I(T_0))$.
 \endproclaim

Since  the logarithmic derivative is not invariant under a M\"obius transformation, (the first assertion of) Theorem 2.2 can not be deduced from Theorem 1.1 directly. We will prove Theorem 2.2 by means of a construction due to Semmes (see [Se1-2]), which is largely different from the approach in our previous paper [Sh2]. Theorem 2.2 only gives a sufficient condition for an increasing  homeomorphism on  the real line  being in  the Weil-Petersson class $\WP(\Bbb R)$. We have shown   in a separated paper (see [STW]) that  this sufficient condition is also a necessary one. Consequently, $\Psi$  is a one-to-one analytic map from the  real Sobolev space $H_{\Bbb R}^{\frac 12}/\Bbb R$ onto the normalized Weil-Petersson class $\WP_0(\Bbb R)$. We will show that the  inverse map $\Psi^{-1}$ is also real analytic.

 \proclaim {Theorem 2.3} $\Psi$  is a one-to-one analytic map from the  real Sobolev space $H_{\Bbb R}^{\frac 12}/\Bbb R$ onto the normalized Weil-Petersson class $\WP_0(\Bbb R)$ whose inverse $\Psi^{-1}$ is also real analytic.
  \endproclaim

  \noindent {\bf Remark 2.1.}\quad Here is an appropriate place to point out why we first prove our main results in the real line case and then come back to the unit circle case. As will be seen later, the main effort of the paper is to prove the real analyticity of the map sending an $H^{\frac 12}$ function to a Weil-Petersson homeomorphism. The proof is based on an important instruction due to Semmes [Se2], which is available on the real line but not on the unit circle. On the other hand, we have a programm to study the Weil-Petersson Teichm\"uller space from several points of view. In the forth-coming work [SW], we will study how the Riemann mapping depends on a  Weil-Petersson quasi-line and  need  the Weil-Petersson  theory on the real line, eg. Theorems 2.2 and 2.3,  parallel to the unit circle case. Theorem 2.3 implies  that the  normalized Weil-Petersson class $\WP_0(\Bbb R)$,   the real model of the  Weil-Petersson Teichm\"uller space $\Cal T$, can be endowed with a real Hilbert manifold structure from $H_{\Bbb R}^{\frac 12}/\Bbb R$   by the correspondence $h\mapsto\log h'$, which is real analytically equivalent to the standard complex Hilbert manifold structure on $\Cal T$   given by Takhtajan-Teo [TT]. This fact  will play an important role in the sequel [SW].

 \head 3 BMO functions
\endhead

In order to prove Theorem 2.2, we need a construction concerning quasiconformal extensions of strongly quasisymmetric homeomorphisms introduced by Semmes [Se1-2], which relies heavily on BMO estimates (see section 4 below). In this section we  recall  some basic definitions and results on BMO functions (see [Gar]).

A locally  integrable function $u\in L^1_{loc}(\Bbb R)$ is said to have bounded mean oscillation and belongs to the space $\BMO$ if
$$\|u\|_{\BMO}\doteq\sup\frac{1}{|I|}\int_I|u(t)-u_I|dt<+\infty,\tag 3.1$$
where the supremum is taken over all finite sub-intervals $I$ of $\Bbb R$, while $u_I$ is the average of $u$ on the interval $I$, namely,
$$u_I=\frac{1}{|I|}\int_Iu(t)dt.\tag 3.2$$ If $u$ also satisfies the condition
$$\lim_{|I|\to 0}\frac{1}{|I|}\int_I|u(t)-u_I|dt=0,$$
we say $u$ has vanishing mean oscillation and belongs to the space $\VMO$. We can define BMO functions and VMO functions on the unit circle in a similar way.  It is well known that  $H^{\frac 12}\subset \VMO$, and the inclusion map is continuous (see [Zh]). In the following, we denote by $\BMO_{\Bbb R}$ the set of all real-valued BMO functions.

We need some basic results on BMO functions. For simplicity, we fix some notations.  $C$,
$C_1$, $C_2$ $\cdots$ will denote universal constants that might
change from one line to another, while $C(\cdot)$, $C_1(\cdot)$,
$C_2(\cdot)$ $\cdots$ will denote constants that depend only on the
elements put in the brackets.
The notation $A \asymp B$ means that there is a positive constant $C$ independent of $A$ and $B$ such that $A/C\le
B\le CA$. The notation $A\lesssim B$ $(A\gtrsim B)$ means that there is a positive constant $C$ independent of $A$ and $B$ such that
$A\le CB$ $ (A\ge CB)$. By the well-known theorem of John-Nirenberg for BMO functions (see [Gar]), there exist two universal positive constants $C_1$ and $C_2$ such that for any BMO function $u$, any subinterval $I$ of $\Bbb R$  and any $\lambda>0$, it holds that
 $$\frac{\left|\{t\in I:|u(t)-u_I|\ge\lambda\}\right|}{|I|}\le C_1\exp\left(\frac{-C_2\lambda}{\|u\|_{\BMO}}\right).\tag 3.3$$
  By  Chebychev's inequality, we obtain that for $u$ with $\|u\|_{\BMO}<C_2$,
 $$
 \align
 \frac{1}{|I|}\int_{I}(e^{|u-u_{I}|}-1)dt&=\frac{1}{|I|}\int_0^{\infty}\left|\{t\in I:|u-u_{I}|\ge\lambda\}\right|d(e^{\lambda}-1) \\
 &\le C_1\int_0^{\infty}e^\lambda\exp\left(\frac{-C_2\lambda}{\|u\|_{\BMO}}\right)d\lambda\tag 3.4\\
 &\le  \frac{C_1\|u\|_{\BMO}}{C_2-\|u\|_{\BMO}}.
 \endalign
 $$
Similarly, for any $p\ge 1$ we have
$$
  \frac{1}{|I|}\int_{I}|u-u_{I}|^pdt\lesssim C(p)\|u\|^p_{\BMO}. \tag 3.5
 $$

 \proclaim {Lemma 3.1} Let $\phi$ be a $C^{\infty}$  function on the real line which is supported on $[-1, 1]$ and satisfies $\int_{\Bbb R}\phi(x)dx=1$. Set $\phi_y(x)=y^{-1}\phi(y^{-1}x)$ for $y>0$, and
 $$\phi_y\ast v(x)=\int_{\Bbb R}\phi_y(x-t)v(t)dt.\tag 3.6$$
 Then it holds that $$|\phi_y\ast e^u|\asymp |e^{\phi_y\ast u}|\tag 3.7$$
 when $\|u\|_{\BMO}$ is small.
 \endproclaim

 \demo{Proof} Lemma 3.1 appears in [Se2]. For completeness and for convenience of later use, we write down  the detailed proof here. Actually, besides Lemma 3.1 itself, the following inequalities (3.8) and (3.9) will also be used in the proof of Theorem 2.2.

 For $x\in\Bbb R$ and $y>0$, consider $I=[x-y, x+y]$ so that
 $$u_I=\frac{1}{2y}\int^{x+y}_{x-y}u(t)dt.$$
Since $\int_{\Bbb R}\phi(x)dx=1$, which implies that $\int_{\Bbb R}\phi_y(x)dx=1$, we obtain
 $$|\phi_y\ast u(x)-u_I|=|\phi_y\ast(u-u_I)(x)|\le C(\phi) \frac{1}{|I|}\int_I|u(t)-u_I|dt\lesssim\|u\|_{\BMO}.\tag 3.8$$
 Since $|e^z-1|\le|ze^z|\le |z|e^{|z|}$, we have
 $$
 \align
 \frac{1}{|I|}\int_I|e^{u(t)-\phi_y\ast u(x)}-1|dt&\le\frac{1}{|I|}\int_I|e^{u(t)-\phi_y\ast u(x)}||u(t)-\phi_y\ast u(x)|dt\\
 &\le\frac{|e^{u_I-\phi_y\ast u(x)}|}{|I|}\int_I|e^{u(t)-u_I}|(|u(t)-u_I|+|u_I-\phi_y\ast u(x)|dt.\\
 \endalign
 $$
 Using H\"older inequality, we conclude from (3.4), (3.5) and (3.8) that
 $$
 \frac{1}{|I|}\int_I|e^{u(t)-\phi_y\ast u(x)}-1|dt\lesssim\|u\|_{\BMO}\tag 3.9$$
 when $\|u\|_{\BMO}$ is small. Noting that
 $$\phi_y\ast e^u(x)-e^{\phi_y\ast u(x)}=e^{\phi_y\ast u(x)}\phi_y\ast(e^{u-\phi_y\ast u(x)}-1)(x),$$
 we obtain
 $$\align
 |\phi_y\ast e^u(x)-e^{\phi_y\ast u(x)}|&=|e^{\phi_y\ast u(x)}||\phi_y\ast(e^{u-\phi_y\ast u(x)}-1)(x)|\\
 &\lesssim \frac{|e^{\phi_y\ast u(x)}|}{|I|}\int_I|e^{u(t)-\phi_y\ast u(x)}-1|dt,
 \endalign$$
 which implies by (3.9)  the required relation (3.7). \quad$\square$
 \enddemo

 \head 4 Semmes' construction revisited
\endhead
 We begin this section with a basic result of Coifman-Meyer [CM]. For $u\in\BMO$ on the real line, set
$$\gamma_u(x)=\frac{\int_0^xe^{u(t)}dt}{\int_0^1e^{u(t)}dt},\quad x\in\Bbb R.\tag 4.1$$
Coifman-Meyer [CM] showed that $\gamma_u$ is a strongly quasisymmetric homeomorphism from the real line $\Bbb R$ onto a chord-arc curve $\Gamma_u=\gamma_u(\Bbb R)$ when $\|u\|_{\BMO}$ is small. If, in addition, $u$ is real-valued, then $\gamma_u$ is a strongly quasisymmetric homeomorphism of $\Bbb R$ onto itself. Recall that   a sense preserving homeomorphism $h$ on $\Bbb R$ is
strongly quasisymmetric if it is locally absolutely continuous so that $|h'|$ belongs to the class of weights $A^{\infty}$ introduced
by Muckenhoupt (see [Gar]) and it maps $\Bbb R$ onto a chord-arc curve passing through the point at infinity (see [Se2]).

In an important paper [Se2], Semmes showed that, when $\|u\|_{\BMO}$ is small,  $\gamma_u$ can be extended to a quasiconformal mapping to the whole plane whose Beltrami coefficient satisfies certain Carleson measure condition. To be precise, let $\varphi$ and $\psi$ be two $C^{\infty}$ real-valued function on the real line  supported on $[-1, 1]$ such that $\varphi$ is even, $\psi$ is odd and $\int_{\Bbb R}\varphi(x)dx=1$, $\int_{\Bbb R}\psi(x)xdx=1$. Define
$$\rho(x, y)=\rho_u(x, y)=\varphi_y\ast \gamma_u(x)-i\psi_y\ast \gamma_u(x),\quad  z=x+iy\in\Bbb U.\tag 4.2$$
Semmes proved that $\rho$ is a quasiconformal mapping from the upper half plane $\Bbb U$ onto the left domain bounded by $\Gamma_u$ when    $\|u\|_{\BMO}$ is small. Furthermore, when $u$ is real-valued, $\rho$ is a quasiconformal mapping of $\Bbb U$ onto itself and is quasi-isometric under the
Poincar\'e metric $|dz|/y$.
 In fact, there exist two $C^{\infty}$ functions
 $\alpha$ and $\beta$   on the real line which are supported on $[-1, 1]$ and satisfy $\int_{\Bbb R}\alpha(x)dx=0$, $\int_{\Bbb R}\beta(x)dx=1$ such that
$$\overline{\partial}\rho(z)=\alpha_y\ast e^u(x),\, {\partial}\rho(z)=\beta_y\ast e^u(x),\quad z=x+iy\in\Bbb U.\tag 4.3$$
It follows from Lemma 3.1 that the Beltrami coefficient $\mu$ of $\rho$ satisfies
$$
\align
|\mu(z)|&=\frac{|\overline{\partial}\rho(z)|}{|{\partial}\rho(z)|}=\frac{|\alpha_y\ast e^u(x)|}{|\beta_y\ast e^u(x)|}\\
&\asymp\frac{|\alpha_y\ast e^u(x)|}{|e^{\beta_y\ast u(x)}|}\\
&=|\alpha_y\ast e^{u-\beta_y\ast u(x)}(x)|\\
&=|\alpha_y\ast (e^{u-\beta_y\ast u(x)}-1)(x)|\tag 4.4\\
&\lesssim\frac{1}{2y}\int^{x+y}_{x-y}|e^{u(t)-\beta_y\ast u(x)}-1)|dt\\
&\lesssim\|u\|_{\BMO}
\endalign
$$
if $\|u\|_{\BMO}$ is small, by (3.9).

 \head 5 Proof of Theorem 2.2 (first part)
\endhead

We first prove the following result.
\proclaim{Lemma 5.1} There exists some universal constant $\delta>0$ such that,  for any  $u\in H^{\frac 12}$ with $\|u\|_{H^{\frac 12}}<\delta$, the mapping $\rho=\rho_u$ defined by (4.2) is quasiconformal whose Beltrami coefficient $\mu$ satisfies $\|\mu\|_{\WP}\lesssim\|u\|_{H^{\frac 12}}$ and thus belongs to the class $\Cal M(\Bbb U)$.
\endproclaim
\demo{Proof} By the continuity of the inclusion $H^{\frac 12}\to\BMO$, we conclude that there exists some universal constant $\delta>0$ such that,  for any  $u\in H^{\frac 12}$ with $\|u\|_{\frac 12}<\delta$, the mapping $\rho=\rho_u$ defined by (4.2) is quasiconformal. It remains to show that $\mu\in\Cal M(\Bbb U)$.

For $z=x+iy\in\Bbb U$, set  $I=[x-y, x+y]$ as before so that
 $$u_I=\frac{1}{2y}\int^{x+y}_{x-y}u(t)dt.$$
Noting that $\int_{\Bbb R}\alpha(x)dx=0$, and  $|e^z-1|\le|ze^z|\le |z|e^{|z|}$, we conclude by (4.3) that
$$
\align
|\overline{\partial}\rho(z)|&=|\alpha_y\ast e^u(x)|=\left|\int_{\Bbb R}\alpha_y(x-t)e^{u(t)}dt\right|\\
&=\left|\int_{\Bbb R}\alpha_y(x-t)(e^{u(t)}-e^{u(x)})dt\right|\\
&=\left|\int_{\Bbb R}\alpha_y(x-t)(e^{u(t)-u(x)}-1)e^{u(x)}dt\right|\\
&\le\int_{\Bbb R}|\alpha_y(x-t)||u(t)-u(x)||e^{u(t)}|dt\\
&\lesssim\frac{1}{|I|}\int_I|u(t)-u(x)||e^{u(t)}|dt.
\endalign
$$
On the other hand, since $\int_{\Bbb R}\beta(x)dx=1$, we conclude by Lemma 3.1 and (4.3) that
$$|{\partial}\rho(z)|=|\beta_y\ast e^u(x)|\asymp|e^{\beta_y\ast u(x)}|.$$
Thus,
$$|\mu(z)|=\frac{|\overline{\partial}\rho(z)|}{|{\partial}\rho(z)|}\lesssim\frac{1}{|I|}\int_I|u(t)-u(x)||e^{u(t)-\beta_y\ast u(x)}|dt.$$
By H\"older inequality, we conclude by (3.9) that
$$\align
|\mu(z)|^2&\lesssim\frac{1}{|I|^2}\int_I|u(t)-u(x)|^2dt\int_I|e^{u(t)-\beta_y\ast u(x)}|^2dt\\
&\lesssim \frac{1}{|I|}\int_I|u(t)-u(x)|^2dt\tag 5.1\\
&\lesssim \frac{1}{y}\int^y_{-y}|u(t+x)-u(x)|^2dt.
\endalign$$
Consequently,
$$
\align
\iint_{\Bbb U}\frac{|\mu(z)|^2}{y^2}dxdy&\lesssim\iint_{\Bbb U}\int^y_{-y}\frac{|u(t+x)-u(x)|^2}{y^3}dtdxdy\\
&=\int^{+\infty}_{-\infty}dx\int^{+\infty}_0\frac{dy}{y^3}\int^y_{-y}|u(t+x)-u(x)|^2dt\\
&=\int^{+\infty}_{-\infty}dx\int^{+\infty}_0\frac{dy}{y^3}\int^y_{0}(|u(x+t)-u(x)|^2+|u(x-t)-u(x)|^2)dt\\
&=\int^{+\infty}_{-\infty}dx\int^{+\infty}_{0}(|u(x+t)-u(x)|^2+|u(x-t)-u(x)|^2)dt\int^{+\infty}_t\frac{dy}{y^3}\tag 5.2\\
&=\int^{+\infty}_{-\infty}dx\int^{+\infty}_{0}\frac{|u(x+t)-u(x)|^2+|u(x-t)-u(x)|^2}{2t^2}dt\\
&=\int^{+\infty}_{-\infty}dx\int^{+\infty}_{-\infty}\frac{|u(x+t)-u(x)|^2}{2t^2}dt\\
&\asymp\|u\|^2_{H^{\frac 12}}.
\endalign
$$\quad$\square$

\enddemo

\proclaim{Corollary 5.1} Let $h$ be an increasing    and  locally absolutely continuous homeomorphism from  the real line  onto itself  such that $\|\log h'\|_{H^{\frac 12}}<\delta$.  Then $h$  can be extended to a quasiconformal mapping to the upper half plane which is quasi-isometric  under the Poincar\'e metric $|dz|/y$ and has Beltrami coefficient in $\Cal M(\Bbb U)$. In particular, $h$ belongs the Weil-Petersson class $\WP(\Bbb R)$.
 \endproclaim

To prove (the first part of) Theorem 2.2, we will decompose a homeomorphism $h$ with finite $\|\log h'\|_{H^{\frac 12}}$ into homeomorphisms $h_j$ with small norms $\|\log h'_j\|_{H^{\frac 12}}$. We need some preliminary results. The first is about the pull-back operator induced by a quasisymmetric homeomorphism. Recall that an  increasing     homeomorphism $h$ from  the real line  onto itself is said to be quasisymmetric if there exists a (least) positive constant $C(h)$, called the quasisymmetric constant of $h$,  such
that $|h(I_1)|\le C(h)|h(I_2)|$
for all pairs of adjacent intervals  $I_1$ and $I_2$ on $\Bbb R$ with the same length $|I_1|=|I_2|$. A strongly quasisymmetric homeomorphism is obviously quasisymmetric. We have the following well-known result.

\proclaim{Proposition 5.1 ([BA], [NS])} Let $h$ be an increasing     homeomorphism $h$ from  the real line  onto itself. Then the pull-back operator $P_h$ defined by $P_hu=u\circ h$ is a bounded operator from $H^{\frac 12}$ into itself if and only if $h$ is quasisymmetric. \endproclaim

\proclaim{Lemma 5.2} Let $h$ be an increasing   and  locally absolutely continuous homeomorphism from  the real line  onto itself  such that $\|\log h'\|_{H^{\frac 12}}<\infty$.  Then  $h$ is strongly quasisymmetric. \endproclaim
\demo{Proof} Consider the Cayley transformation $\gamma(z)=\frac{z-i}{z+i}$ from the upper half plane $\Bbb U$ onto the unit disk $\Delta$. Since $\log h'$ is in $H^{\frac 12}$ on the real line, $\log h'\circ \gamma^{-1}$ is in $H^{\frac 12}$ on the unit circle and consequently in VMO on the unit circle, which implies that $\log h'\circ \gamma^{-1}$ can be approximated by a sequence of bounded functions $(u_n)$ on the unit circle under the BMO norm (see [Gar]). Thus, $\log h'$ can be approximated by the bounded functions $u_n\circ \gamma$ on the real line under the BMO norm. By Lemma 1.4 in [Pa] stating  that an increasing   and  locally absolutely continuous homeomorphism $g$ from  the real line  onto itself is strongly quasisymmetric if $\log g'$ can be approximated by  bounded functions on the real line under the BMO norm, we conclude  that $h$ is strongly quasisymmetric. \quad$\square$
\enddemo

\noindent {\bf Proof of Theorem 2.2 (first part)}\quad Let $h$ be an increasing    and  locally absolutely continuous homeomorphism from  the real line  onto itself  such that $\log h'$ belongs to the Sobolev  class $H^{\frac 12}$. Without loss of generality, we assume $h(0)=0$. For each real number $t\in [0,1]$, set
$$h_t(x)=\int_0^x({h'}(s))^tds,\,x\in\Bbb R.\tag 5.3$$
Then $h_t$ is an increasing    and  locally absolutely continuous homeomorphism from  the real line  onto itself  with  $h_0=\id$, $h_1=h$, and $\log h'_t=t\log h'$, which implies by Lemma 5.2 that $ h_t$ is strongly quasisymmetric. Noting that  for any fixed $t\in [0, 1]$,
$$\|\log(h_{s}\circ h^{-1}_{t})'\|_{H^{\frac 12}}=\|(\log h'_{s}-\log h'_{t})\circ h^{-1}_t\|_{H^{\frac 12}}=|s-t|\|P^{-1}_{h_{t}}\log h'\|_{H^{\frac 12}},$$
we conclude by Proposition 5.1 that there exists a neighbourhood $I_t$ such that $\|\log(h_{s}\circ h^{-1}_{t})'\|_{H^{\frac 12}}<\delta$ when $s\in I_t$. By compactness, we conclude that there exists a sequence of finite numbers $0=t_0<t_1<t_2<\cdots<t_n<t_{n+1}=1$ such that $\|\log(h_{t_j}\circ h^{-1}_{t_{j+1}})'\|_{H^{\frac 12}}<\delta$ for $j=0, 1, 2, \cdots, n-1, n$.  Since $\WP(\Bbb R)$ is a group{\footnote{Cui [Cu] first proved that $\WP(S^1)$,  $\WP(S^1)/\Rot(S^1)$ and $\WP(S^1)/\text{M\"ob}(S^1)$  are all groups (see also [TT]). This can also be seen by means of Theorem 1.1 and Proposition 5.1. Consider  the Cayley transformation $\gamma(z)=\frac{z-i}{z+i}$ from the upper half plane $\Bbb U$ onto the unit disk $\Delta$. Then the correspondence $g\mapsto h=\gamma^{-1}\circ g\circ \gamma$ induces  a one-to-one  from $\WP(S^1)/\Rot(S^1)$ onto $\WP(\Bbb R)$ when $\WP(S^1)/\Rot(S^1)$ is considered as the sub-class of $\WP(S^1)$ of all mappings $h$ with $h(1)=1$. This already implies that $\WP(\Bbb R)$ is also a group. }}, and
$$h^{-1}=(h_{t_0}\circ h^{-1}_{t_1})\circ (h_{t_1}\circ h^{-1}_{t_2})\circ \cdots\circ (h_{t_n}\circ h^{-1}_{t_{n+1}}),$$
We conclude by Corollary 5.1  that  $h\in\WP(\Bbb R)$. \quad$\square$

 \head 6 Proof of Theorem 2.2 (second part) and Theorem 2.3
\endhead

\proclaim{Lemma 6.1} Let $H^{\frac 12}_{\delta}=\{u\in H^{\frac 12}: \|u\|_{H^{\frac 12}}<\delta\}$, where $\delta$ is the universal constant obtained in Lemma 5.1. For $u\in H^{\frac 12}_{\delta}$, let $\Lambda(u)$ denote the Beltrami coefficient for the quasiconformal mapping $\rho_u$ defined by (4.2). Then $\Lambda: H^{\frac 12}_{\delta}\to\Cal M(\Bbb U)$ is holomorphic.
\endproclaim

 \demo{Proof} Since $\Lambda$ is bounded in $H^{\frac 12}_{\delta}$, it is sufficient to show that, for each fixed pair of $(u, v)$ with $u\in H^{\frac 12}_{\delta}$, $v\in H^{\frac 12}$, $\tilde\Lambda(t)\doteq\Lambda(u+tv)$ is holomorphic in a small neighbourhood of $t=0$ in the complex plane. To do so, choose
 $$0<\epsilon<\frac{\delta-\|u\|_{H^{\frac 12}}}{2\|v\|_{H^{\frac 12}}}$$
so that $u+tv\in H^{\frac 12}_{\delta}$ when $|t|\le 2\epsilon$. We conclude by (4.3) that  $\tilde\Lambda(t)(z)$ is holomorphic in $|t|\le 2\epsilon$ for fixed $z\in\Bbb U$. For $|t_0|<\epsilon$, $|t|<\epsilon$, Cauchy formula yields that
$$\align
\left|\frac{\tilde\Lambda(t)(z)-\tilde\Lambda(t_0)(z)}{t-t_0}-\frac{d}{dt}|_{t=t_0}\tilde\Lambda(t)(z)\right|
&=\frac{|t-t_0|}{2\pi }\left|\int_{|\zeta|=2\epsilon}\frac{\tilde\Lambda(\zeta)(z)}{(\zeta-t)(\zeta-t_0)^2}d\zeta\right|\\
&\le\frac{|t-t_0|}{2\pi\epsilon^3}\int_{|\zeta|=2\epsilon}|\tilde\Lambda(\zeta)(z)||d\zeta|.
\endalign
$$
Thus, by (4.4),
$$
\left\|\frac{\tilde\Lambda(t)-\tilde\Lambda(t_0)}{t-t_0}-\frac{d}{dt}|_{t=t_0}\tilde\Lambda(t)\right\|_{\infty}
\le\frac{|t-t_0|}{2\pi\epsilon^3}\int_{|\zeta|=2\epsilon}\|\tilde\Lambda(\zeta)\|_{\infty}|d\zeta|\le C(u, v)|t-t_0|,
$$
and by (5.2),
$$\align
&\iint_{\Bbb U}\frac{1}{y^2}\left|\frac{\tilde\Lambda(t)(z)-\tilde\Lambda(t_0)(z)}{t-t_0}-\frac{d}{dt}|_{t=t_0}\tilde\Lambda(t)(z)\right|^2dxdy\\
&\le\frac{|t-t_0|^2}{4\pi^2\epsilon^6}\iint_{\Bbb U}\frac{1}{y^2}\left(\int_{|\zeta|=2\epsilon}|\tilde\Lambda(\zeta)(z)||d\zeta|\right)^2dxdy\\
&\le\frac{|t-t_0|^2}{\pi\epsilon^5}\iint_{\Bbb U}\int_{|\zeta|=2\epsilon}\frac{|\tilde\Lambda(\zeta)(z)|^2}{y^2}|d\zeta|dxdy\\
&=\frac{|t-t_0|^2}{\pi\epsilon^5}\int_{|\zeta|=2\epsilon}\iint_{\Bbb U}\frac{|\tilde\Lambda(\zeta)(z)|^2}{y^2}dxdy|d\zeta|\\
&\lesssim C(u, v)|t-t_0|^2.
\endalign
$$
 Consequently, the limit
$$\lim_{t\to t_0}\frac{\tilde\Lambda(t)-\tilde\Lambda(t_0)}{t-t_0}=\frac{d}{dt}|_{t=t_0}\tilde\Lambda(t)$$
exists in $\Cal M(\Bbb U)$ and $\Lambda: H^{\frac 12}_{\delta}\to\Cal M(\Bbb U)$ is holomorphic. \quad$\square$
\enddemo

To complete the proof of (the second part of) of  Theorem 2.2, we need to use the  allowable maps introduced in section 2. Let $h_0\in\WP_0(\Bbb R)$ be a normalized mapping in the Weil-Petersson class. Then $g_0=\gamma\circ h_0\circ\gamma^{-1}$ belongs to the Weil-Petersson class $\WP(S^1)$ on the unit circle, where  $\gamma(z)=\frac{z-i}{z+i}$ is the Cayley transformation from the upper half plane $\Bbb U$ onto the unit disk $\Delta$. Cui [Cu] showed that the Douady-Earle [DE] extension $DE(g_0)$ of $g_0$ is a quasiconformal mapping of the unit disk onto itself whose Beltrami coefficient satisfies the condition (1.1). Set $w_0=\gamma^{-1}\circ DE(g_0)\circ\gamma$. Then $w_0$ is a quasiconformal extension of $h_0$ with Beltrami coefficient  $\mu_0\in\Cal
M(\Bbb U)$. Since the Douady-Earle extension $DE(g_0)$ is quasi-isometric under the
Poincar\'e metric $|dz|/(1-|z|^2)$ (see [DE]),  $w_0$ is quasi-isometric under the
Poincar\'e metric $|dz|/y$. Thus $w_0$ induces an allowable map  $w_0^*: {T_0}\to {T_0}$ which is a bi-holomorphic  isomorphism. Now $h_0$ induces a one-to-one mapping $R_{h_0}: \WP_0(\Bbb R)\to \WP_0(\Bbb R)$ which sends $h$ to $h\circ h^{-1}_0$. Clearly,  it holds that $I\circ w^*_0=R_{h_0}\circ I$, which implies that $R_{h_0}: \WP_0(\Bbb R)\to \WP_0(\Bbb R)$ is bi-holomorphic. On the other hand, when $\log h'_0\in H^{\frac 12}_{\Bbb R}$, $h_0$ also induces a bi-holomorphic isomorphism $L_{h_0}: H^{\frac 12}/\Bbb C\to H^{\frac 12}/\Bbb C$ defined by
$$L_{h_0}u=(u-\log h'_0)\circ h^{-1}_0.$$
When restricted to $H^{\frac 12}_{\Bbb R}/\Bbb R$, both $L_{h_0}$ and its inverse are real analytic. Clearly,  $R_{h_0}: \WP_0(\Bbb R)\to \WP_0(\Bbb R)$  and $L_{h_0}: H^{\frac 12}_{\Bbb R}/\Bbb R\to H^{\frac 12}_{\Bbb R}/\Bbb R$ are related by $ R_{h_0}\circ \Psi=\Psi\circ L_{h_0}$. Summarizing these, we obtain
\proclaim{Lemma 6.2} In order to prove that   $\Psi:  H_{\Bbb R}^{\frac 12}/\Bbb R\to \WP_0(\Bbb R)$ is real analytic, it is sufficient to show that  $\Psi$ is real analytic near the base point.
\endproclaim

\noindent {\bf Proof of Theorem 2.2 (second part)}\quad Since $\Psi=I\circ\Phi\circ\Lambda$ on $(H^{\frac 12}_{\delta}\cap H^{\frac 12}_{\Bbb R})/\Bbb R$, we conclude by Lemma 6.1 that $\Psi$ is real analytic near the base point.  Combining this with Lemma 6.2 completes the proof of Theorem 2.2. A direct computation shows that the differential of $\Psi$ at a point $u\in H^{\frac 12}_{\Bbb R}/\Bbb R$  is the linear operator
$$v\mapsto \frac{1}{\left(\int_0^1e^{u(t)}dt\right)^2}\left(\int_0^1e^{u(t)}dt\int_0^xe^{u(t)}v(t)dt-\int_0^1e^{u(t)}v(t)dt\int_0^xe^{u(t)}dt\right)\tag 6.1 $$
for $v\in H^{\frac 12}_{\Bbb R}/\Bbb R$.
\quad$\square$
\vskip 0.2 cm
 \noindent {\bf Proof of Theorem 2.3} As stated in section 2,   $\Psi$  is a one-to-one  map from the  real Sobolev space $H_{\Bbb R}^{\frac 12}/\Bbb R$ onto the normalized Weil-Petersson class $\WP_0(\Bbb R)$. Now  $\Psi$ is real analytic, and its derivative $d_u\Psi: H_{\Bbb R}^{\frac 12}/\Bbb R\to H^{\frac 32}_{\Bbb R}(0, 1)\circ\Psi(u)$  is given by the linear operator (6.1). Here,  $H^{\frac 32}_{\Bbb R}(0, 1)$ is the set of all real-valued $H^{\frac 32}$-functions $ \omega$ with the normalized conditions $\omega(0)=\omega(1)=0$. Recall that  $H^{\frac 32}_{\Bbb R}(0, 1)$ is  the tangent space at the base point of the normalized Weil-Petersson class $\WP_0(\Bbb R)$ (see Theorem 8.1 below), which implies that $H^{\frac 32}_{\Bbb R}(0, 1)\circ\Psi(u)$ is  the tangent space at the  point $\Psi(u)$ of $\WP_0(\Bbb R)$. Clearly, $d_u\Psi: H_{\Bbb R}^{\frac 12}/\Bbb R\to H^{\frac 32}_{\Bbb R}(0, 1)\circ\Psi(u)$ is invertible. In fact, for each $\omega\in H^{\frac 32}_{\Bbb R}(0, 1)\circ\Psi(u)$,
$$(d_u\Psi)^{-1}\omega(x)=\left(\int_0^1e^{u(t)}dt\right)\frac{\omega'(x)}{e^{u(x)}},\quad x\in\Bbb R.\tag 6.2 $$
 Now the invertibility of $d_u\Psi $ implies that the inverse mapping $\Psi^{-1}$ is also real analytic  by the implicit function theorem. \quad$\square$

\vskip 0.2 cm
\noindent {\bf Remark 6.1}  Theorem 2.3 says that, with the standard real Hilbert manifold structure of $H_{\Bbb R}^{\frac 12}/\Bbb R$, $\Psi$  is a one-to-one analytic map from  $H_{\Bbb R}^{\frac 12}/\Bbb R$ onto the normalized Weil-Petersson class $\WP_0(\Bbb R)$ whose inverse $\Psi^{-1}$ is also real analytic. Therefore, there exists a unique complex Hilbert manifold structure on $H_{\Bbb R}^{\frac 12}/\Bbb R$ such that $\Psi$  is a bi-holomorphism  from  $H_{\Bbb R}^{\frac 12}/\Bbb R$ onto $\WP_0(\Bbb R)$. This complex Hilbert manifold structure on $H_{\Bbb R}^{\frac 12}/\Bbb R$ can be assigned as follows. For $u\in H_{\Bbb R}^{\frac 12}/\Bbb R$, define
$$ h_u(x)=\frac{\int_0^xe^{u(t)}dt}{\int_0^1e^{u(t)}dt},\quad x\in\Bbb R$$
as before. By the well-known conformal sewing principle (see [Ah1], [Le], [Na1], [Sh2], [TT]), there exists a pair of quasiconformal mappings $f_u$, $g_u$ on the whole plane $\Bbb C$ which satisfies the following properties:

(1) Both $f_u$ and $g_u$ fix the points $0$, $1$, and $\infty$;

(2) $f_u=g_u\circ h_u$ on the real line;

(3) $f_u$ is conformal in the lower half plane $\Bbb U^*$, with Beltrami coefficient $\mu_1$ in $\Bbb U$ being   square integrable in the Poincar\'e metric, that is,  $\mu_1\in \Cal M(\Bbb U)$;

(4) $g_u$ is conformal in the upper half plane $\Bbb U$, with Beltrami coefficient $\mu_2$ in $\Bbb U^*$ being   square integrable in the Poincar\'e metric, that is,  $\mu_2\in \Cal M(\Bbb U^*)$\footnote {$\Cal M(\Bbb U^*)$ can be defined in the same manner as $\Cal M(\Bbb U)$.}.

\noindent Let
 $\Cal D(\Bbb U^*)$
denote the Dirichlet  space of functions $\phi$ holomorphic in $\Bbb U^*$
with semi-norm
$$
\|\phi\|_{\Cal D(\Bbb U^*)}\doteq\left(\frac{1}{\pi}\iint_{\Bbb U^*}|\phi'(z)|^2dxdy\right)^{\frac 12}.
$$
Then the correspondence $u\mapsto\log f'_u$ induces a one-to-one map from $H_{\Bbb R}^{\frac 12}/\Bbb R$ onto a connected open subset in $\Cal D(\Bbb U^*)/\Bbb C$, which endows $H_{\Bbb R}^{\frac 12}/\Bbb R$ with a complex Hilbert manifold structure. Under this complex Hilbert manifold structure, $\Psi$  is a bi-holomorphism  from  $H_{\Bbb R}^{\frac 12}/\Bbb R$ onto $\WP_0(\Bbb R)$. For more details, see [TT] and also [Sh2].

 \head 7 Proof of Theorems 1.2 and  2.1
\endhead
We first point out the following analogous result to Theorem 1.3. A detailed proof can be found in our paper [HWS].

\proclaim{Theorem 7.1} Under the assumption of Theorem 2.1,   the flow curve $h(t, \cdot)$ of the differential equation $(2.3)$
 satisfies $\log h'(t, \cdot)\in H_{\Bbb R}^{\frac 12}$ and  the mapping $t\mapsto \log h'(t, \cdot)$ from $[0, M]$ into $H_{\Bbb R}^{\frac 12}$ is continuously differentiable such that
$$\frac{d}{dt}\log h'(t, \cdot)=\omega'(t, h(t, \cdot)).\tag 7.1$$
\endproclaim

\noindent {\bf Proof of Theorem 2.1}\quad Since the vector field $\omega(t, \cdot)$ satisfies the normalized condition $\omega(t, 0)=\omega(t, 1)=0$,  the flow curve $h(t, \cdot)$ of the differential equation $(2.4)$ satisfies the condition $h(t, 0)=h(t, 1)-1=0$. Then it holds that $h(t, \cdot)=\Psi(\log h'(t, \cdot)).$  Consequently, $h(t, \cdot)$ is in the normalized Weil-Petersson class $\WP_0(\Bbb R)$ and is continuously differentiable with respect to the Hilbert manifold structure of $\WP_0(\Bbb R)$ by Theorems 2.2 and  7.1. Now since $h(t, \cdot)$ is a smooth curve in the Weil-Petersson class $\WP_0(\Bbb R)$, we have
 $$\left(\frac{d}{dt}h(t, \cdot)\right)(x)=\frac{\partial}{\partial t}(h(t, x)),$$
 which implies  (2.5) from (2.4).
\quad$\square$

\vskip 0.3 cm

\noindent {\bf Proof of Theorem 1.2}\quad Without loss of generality, we may assume that the vector field $\lambda(t, \cdot)$ satisfies the normalized condition $\lambda(t, 1)=\lambda(t, -1)=\lambda(t, -i)=0$ so that  the flow curve $g(t, \cdot)$ of the differential equation $(1.2)$ satisfies the condition $g(t, 1)=1$, $g(t, -1)=-1$, $g(t, -i)=-i$. Consider as above the Cayley transformation $\gamma(z)=\frac{z-i}{z+i}$ from the upper half plane $\Bbb U$ onto the unit disk $\Delta$. Set
$$\omega(t, u)=\frac{\lambda(t, \gamma(u))}{\gamma'(u)},\quad u\in\Bbb R,\tag 7.2$$
and $$h(t, x)=\gamma^{-1}\circ g(t, \gamma(x)),\quad x\in\Bbb R.\tag 7.3$$
By Corollary 8.1 below, we see that  $\omega(t, \cdot)\in C^0([0, M], H_{\Bbb R}^{\frac 32})$ is  a continuous real-valued vector field   on the real line $
\Bbb R$ with $\omega(t, 0)=\omega(t, 1)=0$. A direct computation yields that  $h(t, \cdot)$ is  the flow curve  of the differential equation
$$
\cases
\frac{du}{dt}=\omega(t, u)\\
u(0,x)=x.
\endcases
$$
By Theorem 2.1, $h(t, \cdot)$ is in the normalized Weil-Petersson class $\WP_0(\Bbb R)$ and is continuously differentiable with respect to the Hilbert manifold structure of $\WP_0(\Bbb R)$ such that
$$\frac{d}{dt}h(t, \cdot)=\omega(t, h(t, \cdot)).\tag 7.4$$
Noting that the correspondence $g\mapsto h=\gamma^{-1}\circ g\circ \gamma$ is a real analytic diffeomorphism  from $T_0=\WP(S^1)/\text{M\"ob}(S^1)$ onto $\WP_0(\Bbb R)$, we conclude that $g(t, \cdot)$ is in $T_0=\WP(S^1)/\text{M\"ob}(S^1)$ and is continuously differentiable with respect to the Hilbert manifold structure of $T_0=\WP(S^1)/\text{M\"ob}(S^1)$. Finally, (1.3) follows from (1.2) immediately.
\quad$\square$

 \head 8 Appendix: On the tangent space to $\WP_0(\Bbb R)$
 \endhead
The tangent space to the universal Teichm\"uller space is well understood (see [GL], [GS], [Na2], [Rei]). In this section, we will deal with the tangent space to the Weil-Petersson Teichm\"uller space, showing some results which have been used in previous sections and have independent interests of their own.

Recall that a complex-valued function
$F$ defined in a  domain $\Omega$ is called a quasiconformal
deformation (abbreviated to q.d.)  if it has the generalized
derivative $\bar
\partial F$ such that $\bar \partial F\in L{}^{\infty}(\Omega)$. There
are several reasons for being interested in quasiconformal
deformations because of their close relation with quasiconformal
mappings and Teichm\"uller spaces (see [Ah1], [GL], [GS], [Na1],  [Rei], [WS]) and also of their own interests (see [Ah2], [HSWS], [Re1-2],
[RC], [Sh1], [SLW]). The notion of quasiconformal
deformation is also closely related to the  Zygmund class $\Lambda_*$ in the usual sense (see [Zy]).
Reich-Chen [RC] proved that any Zygmund function $g\in\Lambda_*$ on the unit circle has
a quasiconformal deformation extension to the unit disk and
conversely, any continuous function $g$ on the unit circle which has
a quasiconformal deformation extension to the unit disk must belong
to the Zygmund class $\Lambda_*$ if $g$ also satisfies  the
 condition $\Re\bar {\eta}g(\eta)=0$ for all $\eta\in S^1$.
We will need the following result. A proof of Proposition 8.1 may be founded in our paper [HSWS].

\proclaim{Proposition 8.1} Let $g$ be a continuous function on the unit
circle with the normalized condition $\Re\bar {w}g(w)=0$ on $ S^1$. Then   $g\in H^{\frac 32}$ if and only if  $g$ can be
extended to a quasiconformal deformation $\tilde g$ to the unit disk so that
$$\iint_{\Delta}|\bar\partial\tilde
g(w)|^2(1-|w|^2)^{-2}dudv<+\infty.\tag 8.1$$
\endproclaim

Proposition 8.1 implies that the tangent space to WP$(S^1)$ at the identity consists of precisely the $H^{\frac 32}$ vector fields $\lambda$ on the unit circle (see [NV], [TT] and also [HSWS]), a fact which was already pointed out in section 1. In this section, we will prove the following analogous result.

\proclaim{Theorem 8.1}  The tangent space to $\WP_0(\Bbb R)$ at the identity consists of precisely the $H^{\frac 32}$ real-valued  vector fields $f$ on the real line with the normalized condition $f(0)=f(1)=0$.
\endproclaim

 By the standard Ahlfors-Bers theory of quasiconformal mappings, Theorem 8.1 follows from the following result immediately.

\proclaim{Theorem 8.2} Let $f$ be a real-valued continuous function on the real line.  Then $f\in H^{\frac 32}$ if and only if  $f$ can be
extended to a quasiconformal deformation $\tilde f$ to the upper half plane  so that    $\tilde f(z)=o(z^2)$ as $z\to\infty$ and
$$\iint_{\Bbb U}|\bar\partial\tilde
f(z)|^2y^{-2}dxdy<+\infty,\quad z=x+iy.\tag 8.2$$
\endproclaim

We sketch the standard proof how Theorem 8.1 is deduced from Theorem 8.2 (see [GS], [WS]). Suppose we are given a curve of Weil-Petersson class  mappings
$h^{t}(x)$ ($t>0$ is small) normalized to fix $0$ and $1$, which is the identity
for $t=0$ and differentiable with respect to $t$ for the Hilbert manifold structure on $\WP_0(\Bbb R)=I(T_0)$.
Denote
$$
 h^{t}(x) =x + tf(x) + o(t), \qquad t\rightarrow 0.
$$
Since the natural projection $\Phi$:$ \Cal M(\Bbb U)\to T_0$ is a holomorphic split submersion, we conclude that  there is a differentiable curve  of Beltrami coefficients $\nu_{t}\in\Cal M(\Bbb U)$  such that $h^{t}$ is the
restriction to the real line of the normalized quasiconformal
mapping $f_{\nu_{t}}$. Now there exists some $\mu\in\Cal L(\Bbb U)$ such that
$$\nu_{t} = t\mu + o(t).$$
Consequently,
$$
f_{\nu_{t}}(z) = z + t\dot{f}[\mu](z) + o(t), \qquad t\rightarrow 0. $$
Here  $\dot{f}[\mu]$ satisfies the normalized conditions $\dot{f}[\mu](0)=\dot{f}[\mu](1)=0$, $\dot{f}[\mu](z)=o(z^2)$ as $z\to\infty$,  and is uniquely determined by the condition $\overline{\partial} \dot{f}[\mu]=\mu$ (see [GS]). Noting that $f=\dot{f}[\mu]|_{\Bbb R}$, we conclude by Theorem 8.2 that $f\in H^{\frac 32}$.

Conversely, suppose we are given a function $f\in H^{\frac 32}$
satisfying the normalized condition $f(0)=f(1)=0$. By Theorem 8.2 again,  we deduce that $f$ can
be extended to the upper half plane  to a quasiconformal deformation
$\tilde f$ with $\overline {\partial}-$derivative  $\mu=
\overline {\partial}\tilde f\in\Cal L(\Bbb U)$ and $\tilde f(z)=o(z^2)$ as $z\to\infty$. Set $\mu_t=t\mu$ for small $t>0$. Then
$$
f_{\mu_{t}}(z) = z + t\dot{f}[\mu](z) + o(t), \qquad t\rightarrow 0. $$
Noting that both $\dot{f}[\mu]$ and $\tilde f$ satisfy the normalized conditions $\dot{f}[\mu](0)=\dot{f}[\mu](1)=0$, $\dot{f}[\mu](z)=o(z^2)$ as $z\to\infty$,  and have the same $\overline{\partial}$-derivative $\mu$, we conclude that $\dot{f}[\mu]=\tilde f$. Then,
$$
f_{\mu_{t}}(z) = z + t\tilde f(z) + o(t), \qquad t\rightarrow 0. $$
Set $h^t=f_{\mu_{t}}|_{\Bbb R}$. Then it holds that
$$
 h^{t}(x) =x + tf(x) + o(t), \qquad t\rightarrow 0,
$$
which implies that  $h^t$ is a differentiable curve in $\WP_0(\Bbb R)=I(T_0)$ with the tangent vector $f$.

Before giving the proof of Theorem 8.2, we point out the following corollary, which was already used in the proof of Theorem 1.2.
\proclaim {Corollary 8.1} Let $g$ be a continuous function  on the unit circle with the normalized conditions $g(1)=0$, and $\Re\bar{w}g(w)=0$, and  $\gamma(z)=\frac{z-i}{z+i}$ be the Cayley transformation  from the upper half plane $\Bbb U$ onto the unit disk $\Delta$. Set $f=(g\circ \gamma)/\gamma'$ so that $f$ is a continuous real-valued function on the real line with the normalized condition $f(x)=o(x^2)$ as $x\to\infty$. Then $g\in H^{\frac32}$ on $S^1$ if and only if $f\in H^{\frac32}$ on $\Bbb R$.
\endproclaim
\demo{Proof} For a q.d. extension $\tilde g$ of $g$, $\tilde f=(\tilde g\circ \gamma)/\gamma'$ is a q.d. extension of $f$ with the normalized condition $\tilde f(z)=o(z^2)$ as $z\to\infty$, and vice versa. Moreover, it holds that
$$\bar\partial\tilde f=(\bar\partial\tilde g\circ\gamma)\frac{\overline{\gamma'}}{\gamma'}.$$
Since $\bar\partial\tilde g$ satisfies (8.1) if and only if $\bar\partial\tilde f$ satisfies (8.2), this corollary follows directly from Proposition 8.1 and Theorem 8.2.\quad$\square$
\enddemo

 Now we begin to prove Theorem 8.2. We first recall the following well-known result (see [Zh]).

 \proclaim{Proposition 8.2} Let $\phi$ be analytic in the unit disk. Then it holds that
 $$\iint_{\Delta}|\phi(w)|^2dudv\asymp |\phi(0)|^2+\iint_{\Delta}|\phi'(w)|^2(1-|w|^2)^2dudv.\tag 8.3$$
 \endproclaim
 We show that a similar result also holds on the upper half plane.

 \proclaim{Proposition 8.3} Let $\psi$ be analytic in the upper half plane with $\psi(\infty)=0$. Then it holds that
 $$\iint_{\Bbb U}|\psi(z)|^2dxdy\asymp \iint_{\Bbb U}|\psi'(z)|^2y^2dxdy,\quad z=x+iy.\tag 8.4$$
 \endproclaim

\demo{Proof} Suppose first that $\iint_{\Bbb U}|\psi(z)|^2dxdy<+\infty$. Let $\gamma(z)=\frac{z-i}{z+i}$ be the Cayley transformation  from the upper half plane $\Bbb U$ onto the unit disk $\Delta$ as before. Set $\phi=(\psi\circ{\gamma^{-1}})(\gamma^{-1})'$. Noting that
$$\phi'=(\psi'\circ{\gamma^{-1}}){(\gamma^{-1})'}^2+\psi\circ{\gamma^{-1}}(\gamma^{-1})^{''},$$
we obtain  by Proposition 8.2 that
$$
\align
&\iint_{\Bbb U}|\psi'(z)|^2y^2dxdy\\
&=\iint_{\Delta}|(\psi'\circ{\gamma^{-1}}){(\gamma^{-1})'}^2|^2(1-|w|^2)^2dudv\\
&\lesssim\iint_{\Delta}(|\phi'|^2+|\psi\circ{\gamma^{-1}}(\gamma^{-1})^{''}|^2)(1-|w|^2)^2dudv\\
&\lesssim\iint_{\Delta}(|\phi'(w)|^2(1-|w|^2)^2+|\phi(w)|^2)dudv\\
&\lesssim\iint_{\Delta}|\phi(w)|^2dudv\\
&=\iint_{\Bbb U}|\psi(z)|^2dxdy.
\endalign$$
Here we have used the relation
$$\frac{(\gamma^{-1})^{''}}{(\gamma^{-1})^{'}}(w)=\frac{2}{1-w}.$$

Conversely, suppose that $\iint_{\Bbb U}|\psi'(z)|^2y^2dxdy<+\infty$. Then we have the following reproducing formula (see [GL]):
$$\psi'(z)=\frac{12}{\pi}\iint_{\Bbb U}\frac{v^2\psi'(w)}{(\bar w-z)^4}dudv,\quad w=u+iv,$$
or equivalently,
$$\psi(z)=\frac{4}{\pi}\iint_{\Bbb U}\frac{v^2\psi'(w)}{(\bar w-z)^3}dudv,\quad w=u+iv.$$
Now for any holomorphic function $\varphi$ in the upper half plane with $\iint_{\Bbb U}|\varphi(z)|^2dxdy<+\infty$, we have
$$
\align
&\iint_{\Bbb U}\overline{\psi(z)}\varphi(z)dxdy\\
&=\frac{4}{\pi}\iint_{\Bbb U}\iint_{\Bbb U}\frac{v^2\overline{\psi'(w)}\varphi(z)}{(w-\bar z)^3}dudvdxdy\\
&=-\frac{4}{\pi}\iint_{\Bbb U}v^2\overline{\psi'(w)}dudv\iint_{\Bbb U}\frac{\varphi(z)}{(\bar z-w)^3}dxdy\\
&=-\frac{4}{\pi}\iint_{\Bbb U}v^2\overline{\psi'(w)}dudv\int_{0}^{+\infty}dy\int_{-\infty+iy}^{+\infty+iy}\frac{\varphi(z)}{(z-2iy-w)^3}dz\\
&=-4i\iint_{\Bbb U}v^2\overline{\psi'(w)}dudv\int_0^{+\infty}\varphi^{''}(2iy+w)dy\\
&=2\iint_{\Bbb U}v^2\overline{\psi'(w)}\varphi'(w)dudv,
\endalign
$$
which implies by what we have proved in the first part that
$$
\align
\left|\iint_{\Bbb U}\overline{\psi(z)}\varphi(z)dxdy\right|^2&\le 4\iint_{\Bbb U}|\psi'(z)|^2y^2dxdy\iint_{\Bbb U}|\varphi'(z)|^2y^2dxdy\\
&\lesssim \iint_{\Bbb U}|\psi'(z)|^2y^2dxdy\iint_{\Bbb U}|\varphi(z)|^2dxdy.
\endalign
$$
Consequently,
$$\iint_{\Bbb U}|\psi(z)|^2dxdy\lesssim \iint_{\Bbb U}|\psi'(z)|^2y^2dxdy.$$
\quad$\square$
\enddemo

Now suppose that $f$ is a real-valued continuous function on the
real line, and there exists some constant $\alpha<2$ such that
$f(t)=O(|t|^{\alpha})$ as $t\to\infty$. Following Reich [Re1], set
$$Af(z)=\frac{z^2+1}{i\pi}
\int_{-\infty}^{+\infty}\frac{f(t)}{(t-z)(t^2+1)}dt, \quad  z\in\Bbb U, \tag 8.5$$
and
$$Hf(z)=\frac{(z-\bar z)^3}{2i\pi}
\int_{-\infty}^{+\infty}\frac{f(t)}{(t-z)(t-\bar z)^3}dt, \quad  z\in\Bbb U. \tag 8.6$$
Clearly, $Af$ is analytic on the upper half plane $\Bbb U$, and
$$(Af)^{'''}(z)=\frac{12}{i\pi}
\int_{-\infty}^{+\infty}\frac{f(t)}{(t-z)^4}dt, \quad  z\in\Bbb U. \tag 8.7$$
Reich [Re1] showed that $Hf$ is a $C^{\infty}$ extension of $f$ to $\Bbb U$, and
$$\bar\partial (Hf)(z)=-y^2\overline{(Af)^{'''}(z)}, \quad  z=x+iy\in\Bbb U. \tag 8.8$$
He also showed that  $\Re (Af)$ is continuous extension of $f$ to $\Bbb U$ (see also [SLW]).

\proclaim{Lemma 8.1} Let $f$ be a real-valued continuous function on the
real line with the normalized condition
$f(t)=O(|t|^{\alpha})$ as $t\to\infty$ for  some constant $\alpha<2$. Then $f\in H^{\frac 32}$ if and only if
$$\iint_{\Bbb U}|(Af)^{'''}(z)|^2y^2dxdy,\quad z=x+iy.\tag 8.9$$
\endproclaim
\demo{Proof} Recall that a function $\omega$ on the real line belongs to the class $H^{\frac 12}$ if and only if there exists some harmonic function $\tilde\omega$ on the upper half plane with boundary values $\omega$ and has finite Dirichlet integral $\iint_{\Bbb U}(|\partial\tilde\omega|^2+|\bar\partial\tilde\omega|^2)<+\infty$. Consequently, under the assumption of the lemma, $f\in H^{\frac 32}$ if and only if $\iint_{\Bbb U}|(Af')'-(Af')'(\infty)|^2<+\infty$, which is equivalent to (8.9) by Proposition 8.3 due to the fact that $(Af)^{'''}=(Af')^{''}$. \quad$\square$
\enddemo

 \vskip 0.3 cm

 \noindent {\bf{Proof of Theorem 8.2}}\quad Let $f\in H^{\frac 32}$ be a real-valued continuous function on the real line.  Then we conclude by (8.8) and Lemma 8.1 that $Hf$ is the required  quasiconformal deformation extension of $f$  to the upper half plane. Conversely,  suppose  $f$ can be
extended to a quasiconformal deformation $\tilde f$ to the upper half plane  so that    $\tilde f(z)=o(z^2)$ as $z\to\infty$ and (8.2) holds. Then it holds the following equality (see [Re1] and also [SLW]):
$$\overline{(Af)^{'''}(z)}=-\frac{12}{\pi}\iint_{\Bbb U}\frac{\bar\partial\tilde f(w)}{(w-\bar z)^4}dudv,\quad z=x+iy\in\Bbb U.\tag 8.10$$
A direct computation shows that (8.9) holds by means of (8.2) and (8.10). In fact, by (8.10)        we obtain
$$|{(Af)^{'''}(z)}|^2\le \frac {144}{\pi^2}\iint_{\Bbb U}\frac{|\bar\partial\tilde f(w)|^2}{|w-\bar z|^4}dudv\iint_{\Bbb U}\frac{1}{|w-\bar z|^4}dudv=\frac {36}{\pi y^2}\iint_{\Bbb U}\frac{|\bar\partial\tilde f(w)|^2}{|w-\bar z|^4}dudv,$$
which implies  by (8.2) that
$$\align
\iint_{\Bbb U}|(Af)^{'''}(z)|^2y^2dxdy&\le\frac {36}{\pi}\iint_{\Bbb U}\iint_{\Bbb U}\frac{|\bar\partial\tilde f(w)|^2}{|w-\bar z|^4}dudvdxdy\\
&=\frac {36}{\pi}\iint_{\Bbb U}|\bar\partial\tilde f(w)|^2\iint_{\Bbb U}\frac{1}{|w-\bar z|^4}dxdydudv\\
&=9\iint_{\Bbb U}|\bar\partial\tilde f(w)|^2v^{-2}dudv<+\infty,\, w=u+iv.
\endalign$$
We conclude that $f\in H^{\frac 32}$ by Lemma 8.1 again.  \quad$\square$

\vskip 0.2 cm
\noindent {\bf Acknowledgements} \quad
 The authors would like to thank the  referee for a very careful reading of the manuscript and for several corrections which greatly improves the presentation of the paper.

\enddocument